\documentclass{IEEEtran}
\usepackage{cite}
\usepackage{amsmath,amssymb,amsfonts}
\usepackage{algorithmic}
\usepackage{graphicx}
\usepackage{textcomp}
\usepackage{epstopdf,color}
\allowdisplaybreaks[4]
\newtheorem{theorem}{Theorem}
\newtheorem{lemma}{Lemma}
\newtheorem{remark}{Remark}
\newtheorem{corollary}{Corollary}

\newtheorem{problem}{Problem}
\def\no{\nonumber}
\def\BibTeX{{\rm B\kern-.05em{\sc i\kern-.025em b}\kern-.08em
    T\kern-.1667em\lower.7ex\hbox{E}\kern-.125emX}}
\begin{document}
\title{LQ Optimal Control of First-Order Hyperbolic PDE Systems with Final State Constraints}
\author{Xiaomin~Xue,~
        Juanjuan~Xu,~
        Huanshui~Zhang,~and~Long~Hu
\thanks{This work was supported by the National Natural Science Foundation of China under Grants 61821004, 62250056, 12122110, and the Natural Science Foundation of Shandong Province (ZR2021ZD14, ZR2021JQ24), High-level Talent Team Project of Qingdao West Coast New Area (RCTD-JC-2019-05).}
\thanks{Xiaomin Xue and Juanjuan Xu are with the School of Control Science and Engineering, Shandong University, Jinan, 250061, China (e-mail: xmxue96@163.com, juanjuanxu@sdu.edu.cn).}
\thanks{Huanshui Zhang is with College of Electrical Engineering and Automation, Shandong University of Science and Technology, Qingdao, 266000, China (e-mail: hszhang@sdu.edu.cn).}
\thanks{Long Hu is with School of Mathematics, Shandong University, Jinan, 250100, China (e-mail: hul@sdu.edu.cn).}
}

\maketitle

\begin{abstract}
This paper studies the linear-quadratic (LQ) optimal control problem of a class of systems governed by the first-order hyperbolic partial differential equations (PDEs) with final state constraints.
The main contribution is to present the solvability condition and the corresponding explicit optimal controller by using the Lagrange multiplier method and the technique of solving forward and backward partial differential equations (FBPDEs).
In particular, the result is reduced to the case with zero-valued final state constraints.
Several numerical examples are provided to demonstrate the performance of the designed optimal controller.
\end{abstract}

\begin{IEEEkeywords}
First-order hyperbolic partial differential equations, optimal control, final state constraints, linear quadratic.
\end{IEEEkeywords}

\section{Introduction}
\label{sec:introduction}
Partial differential equations (PDEs) have extensive applications in the real world, which are classified into different types according to evolution natures, consisting of hyperbolic, parabolic, elliptic, etc.
In particular, the first-order hyperbolic PDEs, as one of the most important types of PDEs, are commonly used to characterize various physical phenomena and social behaviors, such as traffic flow \cite{Karafyllis2018feedback}, fuel cell control \cite{Mccain2009dynamic}, cardiovascular flow analysis \cite{Bekiaris20221}.
Therefore, broad attention has been paid to the first-order hyperbolic PDEs and plenty of progress has been made \cite{Coron2021boundary,Meglio2013stabilization,Subbotin2013generalized}.

Two fundamental problems have been stuided for the systems governed by first-order hyperbolic PDEs, one is the optimal control problem, and the other is the controllability problem. On one hand, for the optimal control problem of first-order hyperbolic PDE system, \cite{Aksikas2009lq} and \cite{Moghadam2013boundary} studied the linear-quadratic (LQ) problem and designed the optimal internal control and optimal boundary control based on operator Riccati equations, respectively. In \cite{Xiong2011posteriori}, the optimal control problem was studied by the Galerkin finite element method, which is a kind of numerical methods, and the corresponding posterior error estimates was obtained.
Combining the advantages of the above two methods, \cite{Xue2024linear} provided a computable explicit solution of the LQ optimal control problem.
By reinforcement learning technique, an online policy iteration algorithm for optimal control of first-order hyperbolic PDE systems was designed in \cite{Luo2012online}.
On the other hand, the controllability problem aims to explore whether a dynamical system could be driven by controllers toward the desired final state from an arbitrary initial state.
In \cite{Alabau2017internal} and \cite{Li2003exact}, the exact controllability problems of the first-order hyperbolic PDE systems were studied with the internal control and the boundary control, respectively. The stochastic case was studied in \cite{Lv2014exact,Lu2021mathematical} and the corresponding exact controllability results have been obtained. Based on Carleman estimates, the exact controllability of the time dependent hyperbolic systems was studied in \cite{Klibanov2007exact}.  The minimal time for the exact controllability of first-order hyperbolic systems was presented in \cite{Hu2021minimal,Hu2021null}.

It is noticed that either the optimal control problem without final state constraints, or the controllability problem with the given final state conditions was studied in the above literature.
For the optimal control problem of PDE systems with final state constraints, \cite{Hou2006eigen} studied the parabolic equation, and \cite{Gugat2010penalty} considered the wave equation. However, to the author's best knowledge, there has been no substantial progress on the first-order hyperbolic PDEs.

In this paper, the first-order hyperbolic PDE optimal control problem with final state constraints is studied, where the cost function is considered as quadratic form. The Lagrange multiplier method and the technique of solving forward and backward partial differential equations (FBPDEs) are applied to solve this problem.
The main contributions of this paper are as follows. 1) The solvability condition for the LQ optimal control problem of first-order hyperbolic PDE systems with final state constraints is presented. 2) The corresponding optimal control is given explicitly in the form of state feedback, where a non-homogeneous relationship is established between the optimal costate and the optimal state. Moreover, the result is reduced to a special case of the zero-valued final state constraints.

The rest of the paper is organized as follows. The problem to be addressed is stated in Section \ref{sec2}. The main results in this paper are presented in Section \ref{sec3}. Several numerical examples are provided in Section \ref{sec4}. Finally, section \ref{sec5} concludes the paper.

\textbf{Notations.}
$\mathbb{R}^n$ $(n\geq1)$ denote $n$-dimensional real Euclidean space, $\mathbb{R}=\mathbb{R}^1.$
$x_t(z,t)$ and $x_z(z,t)$ denote the partial derivatives of the variables $t$ and $z$, respectively. For a subset $X$ of $\mathbb{R}^n$ $(n\geq1)$, the following function spaces are useful in the sequel of the paper:
\begin{align*}
C^1(X)=\big\{&x: X\rightarrow \mathbb{R} \big| x\  \mbox{is 1-time continously }\\ & \mbox{differentiable}\big\},
\end{align*}
\begin{align*}
L^2(X)=\big\{&x: X\rightarrow \mathbb{R} \big| x\  \mbox{is Lebesgue measurable,} \\ &\|x\|_{L^2(X)}<\infty\big\},
\end{align*}
where
\begin{align*}
\|x\|_{L^2(X)}:=\left(\int_X \|x\|^2d\sigma\right)^{\frac{1}{2}}.
\end{align*}

\section{Problem Statement}\label{sec2}
Consider the following linear system governed by first-order hyperbolic PDE:
\begin{align}\label{1}
\left\{\aligned
& x_t(z,t)-cx_z(z,t)=ax(z,t)+bu(z,t),\\
&x(z,0)=\varphi(z),\quad
x(l,t)=\phi(t),
\endaligned\right.
\end{align}
for $z\in[0,l], t\in[0,T]$, $l>0, T>0,$ where $x(z,t)\in \mathbb{R}$ is state, $u(z,t)\in \mathbb{R}$ is input control,
$a,$ $b$ are constants, and $c>0$. The initial data  $\varphi(z)\in C^1[0,l]$ and the boundary data $\phi(t)\in C^1[0,T]$ are given with $\varphi(l)=\phi(0)$ for the compatibility at the point $(z,t)=(l,0)$ (i.e., the intersection of the domains of definition).
For any continuously differentiable $u(z,t)$, system (\ref{1}) admits a unique $C^1$ solution, which follows from the PDE well-posedness theory in \cite{Evans2022partial,Li2010controllability}.

The quadratic cost function is considered as
\begin{align}\label{2}
J=&\frac{1}{2}\int_0^T\int_0^l\left[qx^2(z,t)+ru^2(z,t)\right]dzdt\no\\
&+\frac{1}{2}\int_0^lp(z)x^2(z,T)dz,
\end{align}
where constants $q\geq0$, $r>0$, and $p(z)\in C^1[0,l]$ with $p(0)=0$ is a given nonnegative-valued function.

To present the final state constraints on system (\ref{1}) reasonably well,
we introduce the following Riccati equation
\begin{align}\label{26}
\left\{\aligned
& -g_t(z,t)=q+2ag(z,t)+\bar{b}g^2(z,t)-cg_z(z,t),\\
&g(z,T)=p(z),\quad g(0,t)=0,
\endaligned\right.
\end{align}
for $z\in[0,l], t\in[0,T]$, with  $\bar{b}=-\frac{b^2}{r}$. Provided that (\ref{26}) has a solution $g(z,t)$, we define the following functions:
\begin{align}
\bar{a}(z,t)=&a+\bar{b}g(z,t),\label{355}
\end{align}
for $z\in[0,l], t\in[0,T]$, and
\begin{align}
e(z,t)=&e^{-\frac{1}{c}\int_z^l\bar{a}(\xi, t-\frac{\xi-z}{c})d\xi},\label{35}
\end{align}
for $z\in[0,l], t\in[\frac{l-z}{c},T]$.

The addressed LQ optimal control problem with final state constraints is stated as follows.
\begin{problem}\label{P1}
Find the optimal control $u(z,t)$ to minimize (\ref{2}) subject to (\ref{1}) under the final state constraints
\begin{align}\label{3}
x(z,T)=\eta(z),
\end{align}
for $T>\frac{2l}{c}$, $z\in[0,l]$, where $\eta(z)$ is an arbitrary given continuous function with two fixed endpoints:
\begin{align}
\eta(l)=&\phi(T),\\
\eta(0)=&e^{-1}(0,T)\phi(T-\frac{l}{c}).\label{155}
\end{align}
\end{problem}
 \begin{remark}\label{rem1}
The restriction $T>\frac{2l}{c}$ is to ensure that the system can be driven toward the given final state.
 \end{remark}
\begin{remark}\label{rem2}
The requirement $\eta(l)=\phi(T)$ is natural obtained according to (\ref{1}) and (\ref{3}), which is due to the compatibility of $\eta(z)$ and $\phi(t)$ at the point $(x,t)=(l,T)$. Similarly, the requirement of $\eta(0)=e^{-1}(0,T)\phi(T-\frac{l}{c})$ is also to ensure compatibility in one of the solvability conditions of Problem \ref{P1}.
 \end{remark}

\section{Main Results}\label{sec3}
To tackle this type of constrained optimal control problem, the first step is to convert it into a problem without final state constraints by using the Lagrange multiplier method. Specifically, by introducing an augmented functional
\begin{align}\label{4}
\bar{J}=&\frac{1}{2}\int_0^T\int_0^l\left[qx^2(z,t)+ru^2(z,t)\right]dzdt+\frac{1}{2}\int_0^lp(z)\no\\&\times
x^2(z,T)dz
+\int_0^l\gamma(z)[x(z,T)-\eta(z)]dz,
\end{align}
for $T>\frac{2l}{c},$ where $\gamma(z)$ can be an arbitrary constant-valued function, then Problem \ref{P1} can be restated as below:
\begin{problem}\label{P2}
Find the optimal control $u(z,t)$ to minimize (\ref{4}) subject to (\ref{1}).
\end{problem}

Given that Problem \ref{P1} has been converted into Problem \ref{P2}, we first present the solvability result of Problem \ref{P2} as follows.
\begin{lemma}\label{lem1}
For an arbitrary given $\gamma(z)$ in (\ref{4}), if the following FBPDEs are solvable:
\begin{align}
& x_t(z,t)=cx_z(z,t)+ax(z,t)+bu(z,t),\label{15}\\
&-\lambda_t(z,t)=-c\lambda_z(z,t)+a\lambda(z,t)+qx(z,t),\label{16}\\
&x(z,0)=\varphi(z),\quad x(l,t)=\phi(t),\label{17}\\
&\lambda(z,T)=p(z)x(z,T)+\gamma(z),\quad \lambda(0,t)=0,\label{18}\\
&0=ru(z,t)+b\lambda(z,t),\label{19}
\end{align}
then Problem \ref{P2} is solvable with the given $\gamma(z)$, and the optimal controller satisfies (\ref{19}).
\end{lemma}
\textbf{Proof.}
By defining $J(u)$ derived by arbitrary controller $u(z,t)$ with the associated state $x(z,t)$ satisfying (\ref{1}), and defining $J(u^*)$ derived by (\ref{2}) with the controller $u^*(z,t)$ and the state $x^*(z,t)$ satisfying FBPDEs (\ref{15})-(\ref{19}), now we show that
\begin{align}\label{J}
J(u^*)\leq J(u).
\end{align}

In view of cost function (\ref{2}), it follows that
\begin{align}\label{l1}
&2[J(u)-J(u^*)]\no\\
=&\int_0^T\int_0^l\big\{q [x(z,t)-x^*(z,t)]^2+2q [x(z,t)-x^*(z,t)] \no\\&\times x^*(z,t)+r[u(z,t)-u^*(z,t)]^2
+2 r[u(z,t)\no\\&-u^*(z,t)] u^*(z,t)\big\}dzdt
+\int_0^l\big\{p(z)[x(z,T)\no\\&-x^*(z,T)]^2+2p(z)[x(z,T)-x^*(z,T)]x^*(z,T)\big\}dz\no\\&
+2 \int_0^l\gamma(z)[x(z,T)-x^*(z,T)]dz\no\\
\geq
&2\int_0^T\int_0^l\big\{[x(z,t)-x^*(z,t)]qx^*(z,t)\no\\&+[u(z,t)-u^*(z,t)]r u^*(z,t)]dzdt\no\\&
+2\int_0^l\big\{[x(z,T)-x^*(z,T)] p(z)x^*(z,T)\no\\&
+[x(z,T)-x^*(z,T)]\gamma(z)\big\}dz,
\end{align}
where the last inequality is due to the fact that the values of $q$, $r$, and $p(z)$ are nonnegative.
By substituting (\ref{16}), (\ref{18}) and (\ref{19}) into (\ref{l1}), and using integration by parts,  we derive that
\begin{align}\label{l2}
&J(u)-J(u^*)\no\\
\geq
&\int_0^T\int_0^l
\big\{[x(z,t)-x^*(z,t)]_t-c[x(z,t)-x^*(z,t)]_z\no\\&
-a[x(z,t)-x^*(z,t)]
-b[u(z,t)-u^*(z,t)]\big\}\no\\&\times \lambda(z,t)dzdt
+\int_0^l\big\{[x(z,T)-x^*(z,T)]\no\\&\times [\lambda(z,T) -\gamma(z)]
+[x(z,T)-x^*(z,T)]\gamma(z)\big\}dz\no\\
&-\int_0^l\big\{[x(z,t)-x^*(z,t)] \lambda(z,t)\big\}dz\bigg|_{t=0}^{t=T}\no\\&
+\int_0^T\big\{[x(z,t)-x^*(z,t)] c\lambda(z,t)\big\}dt\bigg|_{z=0}^{z=l}.
\end{align}
Combining with (\ref{1}), (\ref{15}), (\ref{17}), and (\ref{18}), it is obtained from (\ref{l2}) that
\begin{align*}
J(u)-J(u^*)\geq0.
\end{align*}
Thus, (\ref{J}) holds. This implies that the controller $u^*(z,t)$ (\ref{19}) is optimal, and thus Problem \ref{P2} is solvable. The proof is now completed.
\hfill$\Box$

Observed that the solvability of Problem \ref{P2} with the arbitrary given $\gamma(z)$ has been provided in Lemma \ref{lem1}. To obtain the solvability of Problem \ref{P1}, in what follows, we define the function
\begin{align}
\gamma(z)=&[e(z,T)\eta(z)-\phi(T-\frac{l-z}{c})]e(z,T)\no\\\
&\times c\bigg[\int_z^l\bar{b}e^2(\xi,T-\frac{\xi-z}{c})d\xi\bigg]^{-1},\label{366}
\end{align}
for $z\in[0,l], t\in[\frac{l-z}{c},T]$,
and further introduce the equation
\begin{align}\label{29}
\left\{\aligned
& -\psi_t(z,t)=\bar{a}(z,t)\psi(z,t)-c\psi_z(z,t),\\
&\psi(z,T)=\gamma(z),\quad \psi(0,t)=0,
\endaligned\right.
\end{align}
for $z\in[0,l], t\in[0,T]$.
Then, the main results of Problem \ref{1} are presented as follows.
\begin{theorem}\label{Thm1}
If (\ref{26}) has a solution $g(z,t)$, (\ref{29}) has a solution $\psi(z,t)$ with the final data $\gamma(z)$ defined by (\ref{366}),
then Problem \ref{P1} is solvable.
With this conditions holding, the optimal control is given by
\begin{align}\label{32}
u^*(z,t)=-\frac{b}{r}[g(z,t)x^*(z,t)+\psi(z,t)],
\end{align}
where the optimal state $x^*(z,t)$ satisfies
\begin{align}\label{33}
\left\{\aligned
&x^*_t(z,t)-c x^*_z(z,t)=\bar{a}(z,t)x^*(z,t)+\bar{b}\psi(z,t),\\
&x^*(z,0)=\varphi(z),\quad
x^*(l,t)=\phi(t),
\endaligned\right.
\end{align}
for $z\in[0,l], t\in[0,T]$.
The corresponding optimal cost is
\begin{align}\label{66}
J^*=&\frac{1}{2}\bigg\{\int_0^T\{c\phi(t)[g(l,t)\phi(t)+\psi(l,t)\}dt\no\\
&+\int_0^l\{\varphi(z)[g(z,0)\varphi(z)+\psi(z,0)\}dz\no\\
&-\int_0^l\eta(z)\gamma(z)dz\bigg\}.
\end{align}
Moreover, the optimal costate and state have the relationship
\begin{align}\label{31}
\lambda^*(z,t)=g(z,t)x^*(z,t)+\psi(z,t).
\end{align}
\end{theorem}
\textbf{Proof.}
See Appendix A for more details of the proof.
\hfill$\Box$
 \begin{remark}\label{rem3}
As explained in Remark \ref{rem2}, the requirement (\ref{155}) of $\eta(0)$ is presented to ensure the compatibility of $\gamma(0)=0$ in equation (\ref{29}), which is necessary for the solvability of equation (\ref{29}).
  \end{remark}
 \begin{remark}\label{rem4}
The solvability of equations (\ref{26}) and (\ref{29}) have been well studied in \cite{Li2010controllability}. For example, one of the sufficient conditions to ensure solvability can be found in Theorem 2.5 in \cite{Li2010controllability}. That is, if the final data $p(z)$ and $\gamma(z)$ of (\ref{26}) and (\ref{29}) are suitably small by choosing $q$, $r$, and $p(z)$ in cost function (\ref{2}), then (\ref{26}) and (\ref{29}) admit unique $C^1$ solutions $g(z,t)$ and $\psi(z,t)$, respectively, for $x\in[0,L]$, $t\in[0,T]$.
  \end{remark}

In Theorem \ref{Thm1}, the optimal control result under the final state data $\eta(z)$ has been provided. Now we reduce the result into the case of zero-valued final state constraints to show the wider application of Theorem \ref{Thm1}.
Specifically, for a class of commonly applied boundary conditions $\phi(t)=0$ (eg. \cite{Tatsien2009exact}), system (\ref{1}) is rewritten as
\begin{align}\label{1111}
\left\{\aligned
& x_t(z,t)-cx_z(z,t)=ax(z,t)+bu(z,t),\\
&x(z,0)=\varphi(z),\quad
x(l,t)=0.
\endaligned\right.
\end{align}
Then the corresponding result to system (\ref{1111}) and the explicit optimal controller are provided in the following corollary.
\begin{corollary}\label{Cor1}
If (\ref{26}) and (\ref{29}) with $\gamma(z)=0$ have solutions $g(z,t)$ and $\tilde{\psi}(z,t)$, then the controller is given by
\begin{align}\label{322}
u^*(z,t)=-\frac{b}{r}[g(z,t)\tilde{x}^*(z,t)+\tilde{\psi}(z,t)],
\end{align}
where the optimal state $\tilde{x}^*(z,t)$ satisfies
\begin{align}\label{333}
\left\{\aligned
&\tilde{x}^*_t(z,t)-c \tilde{x}^*_z(z,t)=\bar{a}(z,t)\tilde{x}^*(z,t)+\bar{b}\tilde{\psi}(z,t),\\
&\tilde{x}^*(z,0)=\varphi(z),\quad
\tilde{x}^*(l,t)=0,
\endaligned\right.
\end{align}
for $z\in[0,l], t\in[0,T]$.
In this case, controller (\ref{322}) is optimal to minimize (\ref{2}), subject to system (\ref{1111}),
and enables that system (\ref{1111}) satisfies the final conditions $x(z,T)=0$.
Furthermore, the corresponding optimal cost is equal to
\begin{align}\label{666}
J^*=\frac{1}{2}\int_0^l\{\varphi(z)[g(z,0)\varphi(z)+\psi(z,0)\}dz.
\end{align}
\end{corollary}
\textbf{Proof.}
By letting $\phi(t)=0$ in system (\ref{1}) and $\gamma(z)=0$ in (\ref{366}), Corollary \ref{Cor1} follows from Theorem \ref{Thm1}.
\hfill$\Box$
\begin{remark}\label{rem5}
According to Corollary \ref{Cor1}, system (\ref{1111}) can be stabilized in finite time by the following extended optimal controller
\begin{align}\label{3222}
u^*(z,t)=\left\{\aligned
&-\frac{b}{r}[g(z,t)\tilde{x}^*(z,t)+\tilde{\psi}(z,t)], \ t\in[0,T],\\
&0,  \ t\in(T, +\infty),
\endaligned\right.
\end{align}
for $z\in[0,l]$, that is, the corresponding closed-loop system satisfies $x(z,t)=0$ for all $t\geq T.$
\end{remark}

\section{Numerical Examples}\label{sec4}
\begin{figure}[htbp]

     \begin{minipage}[t]{0.48\textwidth}
    \centering
     \includegraphics[width=7cm]{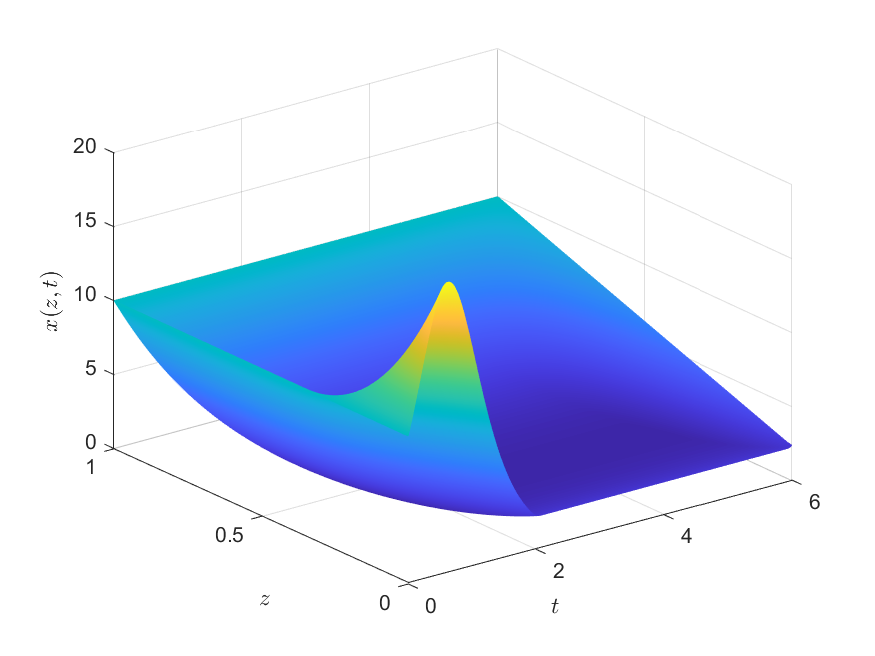}
     \caption{Simulations for the state trajectory $x(z,t)$ with the optimal controller}\label{Fig1}
     \end{minipage}

       \begin{minipage}[t]{0.48\textwidth}
    \centering
     \includegraphics[width=7cm]{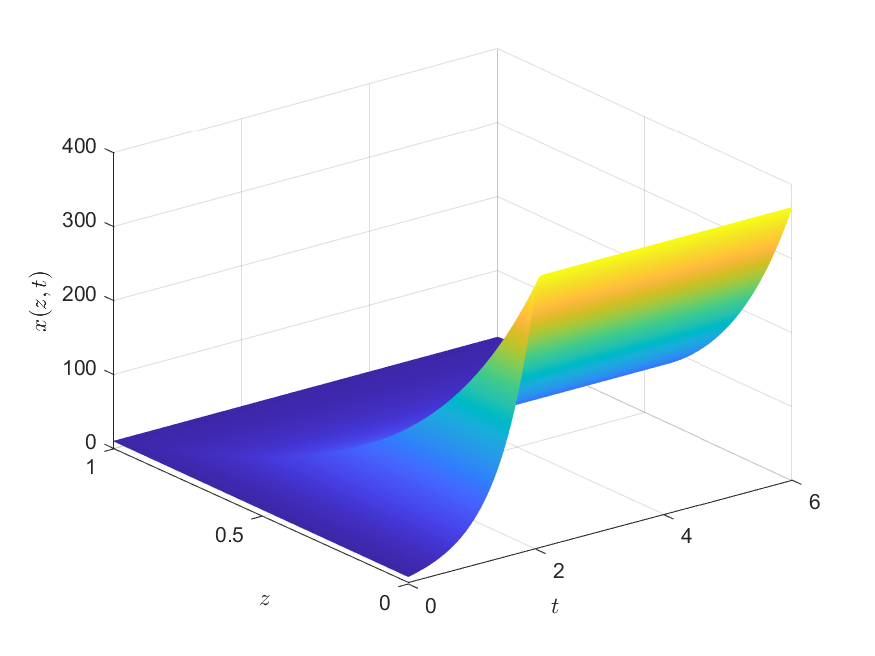}
     \caption{Simulations for the state trajectory $x(z,t)$ without the optimal controller}\label{Fig2}
          \end{minipage}

\end{figure}
\begin{figure}[htbp]

     \begin{minipage}[t]{0.48\textwidth}
    \centering
     \includegraphics[width=7cm]{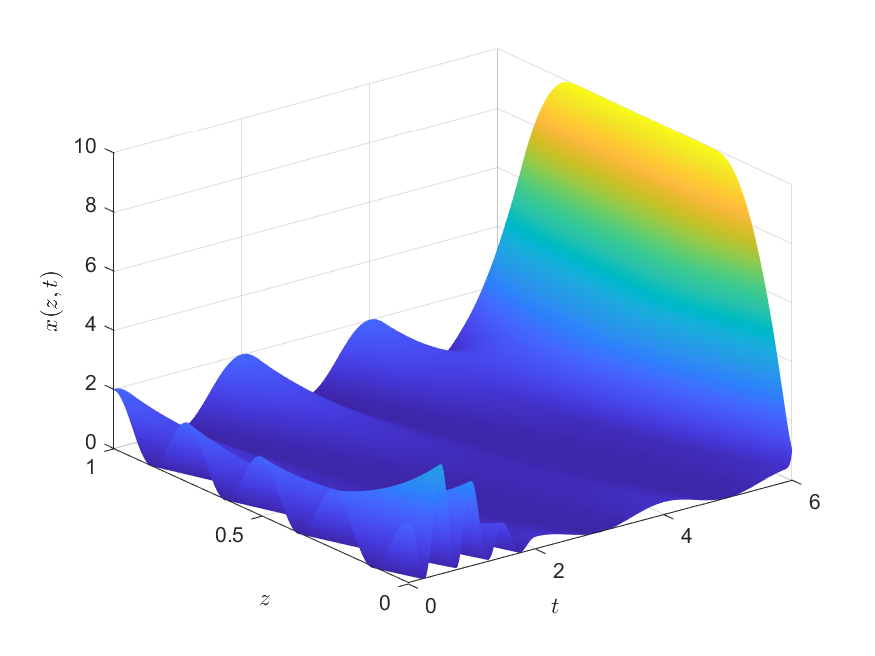}
     \caption{Simulations for the state trajectory $x(z,t)$ with the optimal controller}\label{Fig3}
     \end{minipage}

       \begin{minipage}[t]{0.48\textwidth}
    \centering
     \includegraphics[width=7cm]{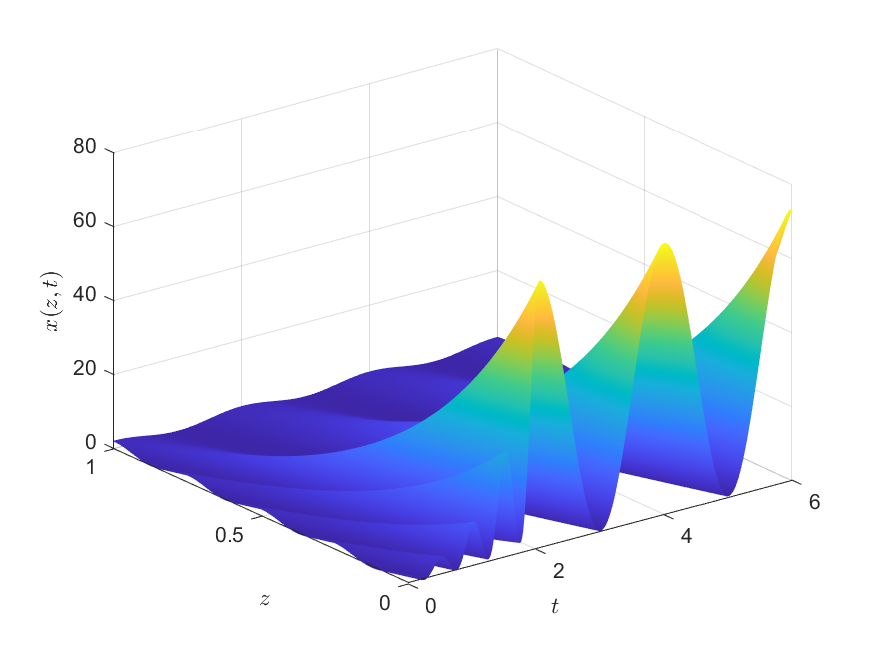}
     \caption{Simulations for the state trajectory $x(z,t)$ without the optimal controller}\label{Fig4}
          \end{minipage}

\end{figure}

In this section, several numerical examples are simulated to verify Theorem \ref{Thm1},
where the simulations are performed on the time interval $[0, 6]$ and the spatial interval $[0, 1]$ through the finite difference method with the step size
$\tau=0.001$ and $ h=0.001.$

Firstly, system (\ref{1}) and cost function (\ref{2}) are considered with the following parameters
$$a=1.8, \ b=1.3, \  c=0.5, \
 q=2, \ r=3, \ p(z)=z.$$
The initial data and boundary data are given by
\begin{align}\label{m1}
\varphi(z)=10,\quad \phi(t)=10,
\end{align}
for $z\in[0, 1], t\in[0, 6]$. The final conditions are considered as
\begin{align}\label{m2}
\eta(z)=\eta(0)+[10-\eta(0)]z, \quad z\in[0, 1].
\end{align}
By solving Riccati equation ({26}) and calculating (\ref{155}), we obtain that $\eta(0)=2.2737$ in Problem \ref{1}.Then, we design the optimal controller by (\ref{32}) and simulate the corresponding state trajectory by (\ref{1}). As shown in Fig. \ref{Fig1},
it can be observed that the state trajectory under the designed controller reaches the given final conditions (\ref{m2}) at the final time $t = 6.$ For comparison, the corresponding state trajectory without controller is depicted in Fig. \ref{Fig2}.

To illustrate the powerful applicability of the result in this paper for various types of data, we present another numerical example by reselecting the initial data and boundary data as
\begin{align}\label{m3}
\varphi(z)=1+cos(8\pi z), \quad\phi(t)=1+cos(\pi t),
\end{align}
and the final condition as the following piecewise function
\begin{align}\label{m4}
\eta(z)=\left\{\aligned
10-16(10-\eta(0))(z-\frac{1}{4})^2, \quad z\in[0,\frac{1}{4}],\\
10, \quad z\in(\frac{1}{4}, \frac{3}{4}],\\
10-128(z-\frac{3}{4})^2, \quad z\in(\frac{3}{4}, 1],
\endaligned\right.
\end{align}
where $\eta(0)=0.4547$ can be obtained by (\ref{26}) and (\ref{155}). Then, by designing the optimal controller according to (\ref{32}),
the state trajectories with and without the designed controller are shown in Fig. \ref{Fig3} and Fig. \ref{Fig4}, respectively. It can be seen that the final state satisfies (\ref{m4}) under the designed controller as shown in Fig. \ref{Fig3}.

\section{Conclusion}\label{sec5}
In this paper, the constrained LQ optimal control problem of first-order hyperbolic PDE systems has been studied. The optimal controllers and the corresponding optimal cost of this problem have been explicitly given. The result has been reduced to the case of zero-valued final state constraints.
Finally, several numerical examples have been simulated to verify the main results.


{\appendices
\section{Proof of Theorem \ref{Thm1}}
\renewcommand{\theequation}{A.\arabic{equation}}

The proof is divided into three steps: (i) Show that (\ref{32}), (\ref{33}), and (\ref{31}) are optimal control, optimal state, and optimal costate, respectively; (ii) Verify the  optimal state (\ref{33}) satisfies final state constraints $x^*(z,T)=\eta(z)$; (iii) Calculate the corresponding optimal cost. The details are as follows.

(i) According to Lemma \ref{lem1}, we only need to show that ($u^*, x^*, \lambda^*$) given by (\ref{32}), (\ref{33}), and (\ref{31}) are the solutions of FBPDEs (\ref{15})-(\ref{19}), in which $g(z,t)$ and $\psi(z,t)$ satisfy equation (\ref{26}) and (\ref{29}) with the definition (\ref{366}).

Firstly, in view of (\ref{32}) and (\ref{31}), it is obviously that $u^*$ defined by (\ref{32}) satisfies (\ref{19}). Secondly, by the definition (\ref{355}) and substituting (\ref{32}) into (\ref{33}), we have
\begin{align}\label{30}
\left\{\aligned
&x^*_t(z,t)-c x^*_z(z,t)=ax^*(z,t)\\
&\qquad\qquad-\frac{b^2}{r}[g(z,t) x^*(z,t)+\psi(z,t)],\\
&x^*(z,0)=\varphi(z),\quad
x^*(l,t)=\phi(t),
\endaligned\right.
\end{align}
which gives that $x^*$ defined by (\ref{33}) satisfies (\ref{15}) and (\ref{17}). Thirdly, by taking partial derivatives with respect to (\ref{31}), it is obtained that
\begin{align}
&-\lambda^*_t(z,t)+c\lambda^*_z(z,t)\no\\
=&-[g(z,t)x^*(z,t)+\psi(z,t)]_t+c[g(z,t)x^*(z,t)\no\\&+\psi(z,t)]_z\no\\
=&[-g_t(z,t)+cg_z(z,t)]x^*(z,t)-g(z,t)[x^*_t(z,t)\no\\&-cx^*_z(z,t)]
-\psi_t(z,t)+c\psi_z(z,t),
\end{align}
which can be calculated by (\ref{30}), (\ref{26}), and (\ref{29}) that
\begin{align}
&-\lambda^*_t(z,t)+c\lambda^*_z(z,t)\no\\
=&[q+2ag(z,t)-\frac{b^2}{r}g^2(z,t)]x^*(z,t)-g(z,t)\no\\
&\times\{ax^*(z,t)-\frac{b^2}{r}[g(z,t)x^*(z,t)+\psi(z,t)]\}\no\\
&+a\psi(z,t)-\frac{b^2}{r}g(z,t)\psi(z,t)\no\\
=&qx^*(z,t)+a\lambda^*(z,t).
\end{align}
This implies that $\lambda^*$ defined by (\ref{31}) satisfies (\ref{16}). Moreover, the final conditions and boundary conditions in (\ref{26}) and (\ref{29}) yield that $\lambda^*$ satisfies (\ref{18}). Therefore, we have shown that ($u^*, x^*, \lambda^*$) satisfies FBPDEs (\ref{15})-(\ref{19}).

(ii) The second step is to verify $x^*(z,T)=\eta(z)$.

By using the invertible transformation $w(z,t)=e(z,t)x^*(z,t)$ for $z\in[0,l], t\in[\frac{l-z}{c},T]$, equation (\ref{33}) of $x(z,t)$ can be reformulated as
\begin{align}\label{37}
&w_t(z,t)-c w_z(z,t)\no\\
=&-\frac{1}{c}\int_z^l\bar{a}_t(\xi, t-\frac{\xi-z}{c})d\xi  e(z,t) x^*(z,t)
+e(z,t) \no\\&\times x^*_t(z,t)-c\bigg\{(-\frac{1}{c})[-\bar{a}(z,t)
+\int_z^l\frac{1}{c}\bar{a}_t(\xi, t-\frac{\xi-z}{c})d\xi]\no\\&\times e(z,t) x^*(z,t)
+e(z,t) x^*_z(z,t)\bigg\}\no\\
=&e(z,t)\bar{b}\psi(z,t),
\end{align}
with $w(l,t)=\phi(t)$.
 By the method of characteristics in \cite{Coron2021boundary}, the solution of (\ref{37}) is written explicitly as
\begin{align}\label{44}
w(z,t)=
&\phi(t-\frac{l-z}{c})+\frac{1}{c}\int_z^l e(\xi,t-\frac{\xi-z}{c})\no\\
&\times\bar{b}\psi(\xi,t-\frac{\xi-z}{c})d\xi,
\end{align}
for $z\in[0,l], t\in[\frac{l-z}{c},T]$.

Noting that $\psi$ in (\ref{44}) is not the final desired result, so we need to further simplify (\ref{44}).
Specifically, according to (\ref{355}), the equation (\ref{29}) of $\psi(z,t)$ is equivalent to
\begin{align}\label{46}
\left\{\aligned
&\psi_t(z,t)-c\psi_z(z,t)+\bar{a}(z,t)\psi(z,t)=0,\\
&\psi(z,T)=\gamma(z),\quad \psi(0,t)=0,
\endaligned\right.
\end{align}
which can be reformulated as
\begin{align}\label{49}
\left\{\aligned
&v_t(z,t)-cv_z(z,t)=0,\\
&v(z,T)=e^{-1}(z,T)\gamma(z),\quad v(0,t)=0,
\endaligned\right.
\end{align}
by the transformation $v(z,t)=e^{-1}(z,t)\psi(z,t)$.
Similar to (\ref{44}), the solution of (\ref{49}) is also obtained by the method of characteristics explicitly as
\begin{align}\label{58}
v(z,t)=e^{-1}(z-c(T-t),T)\gamma(z-c(T-t)),
\end{align}
for $z\in[0,l], t\in[T-\frac{z}{c},T]$. In view of $T>\frac{2l}{c},$ and by substituting $\psi(z,t)=e(z,t)v(z,t)$ with (\ref{58}) into (\ref{44}), we obtain that  $w(z,t)$ can be further expressed as
 \begin{align}\label{62}
w(z,t)=
&\phi(t-\frac{l-z}{c})+\frac{1}{c}\int_z^le^2(\xi,t-\frac{\xi-z}{c})d\xi\no\\
&\times \bar{b}e^{-1}(z-c(T-t),T)\gamma(z-c(T-t)),
\end{align}
for $z\in[0,l], t\in[\frac{2l}{c},T]$.

Therefore, from (\ref{62}) and the transformation $ x^*(z,t)=e^{-1}(z,t)w(z,t)$, we have
 \begin{align}\label{622}
x^*(z,t)=&e^{-1}(z,t)\phi(t-\frac{l-z}{c})\no\\
& +e^{-1}(z,t)\int_z^le^2(\xi,t-\frac{\xi-z}{c})d\xi \no\\
&\times\frac{\bar{b}}{c}e^{-1}(z-c(T-t),T)\gamma(z-c(T-t)),
\end{align}
for $z\in[0,l], t\in[\frac{2l}{c},T]$.

When $t=T$, it is observed from (\ref{622}) that $x^*(z,T)$ is equal to
 \begin{align}\label{64}
x^*(z,T)=
&e^{-1}(z,T)\phi(T-\frac{l-z}{c})+[e^{-1}(z,T)]^2\gamma(z)\no\\
&\times\frac{\bar{b}}{c}\int_z^le^2(\xi,T-\frac{\xi-z}{c})d\xi,
\end{align}
for all $z\in[0,l].$
This implies that $x(z,T)=\eta(z)$ by the definition (\ref{366}) of $\gamma(z)$.

(iii) The third step is to calculate the corresponding optimal cost.

Define the Lyapunov function
\begin{align}\label{71}
L(t)=\int_0^l x^*(z,t)[g(z,t)x^*(z,t)+\psi(z,t)]dz.
\end{align}
On one hand, it is obtained from (\ref{71}) that
\begin{align}\label{72}
&\frac{d}{dt}\int_0^T L(t)dt=L(T)-L(0)\no\\
=&\int_0^l x^*(z,T)[g(z,T)x^*(z,T)+\psi(z,T)]dz\no\\&-\int_0^l x^*(z,0)[g(z,0)x^*(z,0)+\psi(z,0)]dz\no\\
=&\int_0^l x^*(z,T)[p(z)x^*(z,T)+\gamma(z)]dz\no\\&-\int_0^l \varphi(z)[g(z,0)\varphi(z)+\psi(z,0)]dz,
\end{align}
where the last equality holds because of the final conditions in (\ref{26})-(\ref{29}) and initial conditions in (\ref{30}).
On the other hand, by using  (\ref{30}), (\ref{26}), (\ref{29}), and (\ref{32}), we have
\begin{align}\label{73}
&\frac{d}{dt}\int_0^T L(t)dt\no\\
=&\int_0^T\int_0^lx^*_t(z,t)[2g(z,t)x^*(z,t)+\psi(z,t)]dzdt\no\\&+\int_0^T\int_0^l x^*(z,t)[g_t(z,t)x^*(z,t)+\psi_t(z,t)]dzdt\no\\
=&\int_0^T\int_0^l \{cx^*_z(z,t)+ax^*(z,t)-\frac{b^2}{r}[g(z,t)x^*(z,t)\no\\&+\psi(z,t)]\} [2g(z,t)x^*(z,t)+\psi(z,t)]dzdt\no\\&
-\int_0^T\int_0^l x^*(z,t)\{[q+2ag(z,t)-\frac{b^2}{r}g^2(z,t)\no\\&-cg_z(z,t)]x^*(z,t)+[a\psi(z,t)-\frac{b^2}{r}g(z,t)\psi(z,t)\no\\&-c\psi_z(z,t)]\}dzdt\no\\
=&\int_0^T\int_0^l\{cx^*_z(z,t)\psi(z,t)+c x^*(z,t)\psi_z(z,t)\no\\&
+ cx^*_z(z,t)2g(z,t)x^*(z,t)+ cx^*(z,t)g_z(z,t)x^*(z,t)\no\\&
- x^*(z,t)qx^*(z,t)
- u^*(z,t)ru^*(z,t)\}
dzdt,
\end{align}
which is further calculated according to integration by parts and the boundary conditions in (\ref{30}), (\ref{26}), and (\ref{29}) as
\begin{align}\label{74}
&\frac{d}{dt}\int_0^T L(t)dt\no\\
=&\int_0^T[cx^*(l,t)\psi(l,t)-cx^*(0,t)\psi(0,t)]dt\no\\&
+\int_0^T[ cx^*(l,t)g(l,t)x^*(l,t)-cx^*(0,t)g(0,t)x^*(0,t)]dt\no\\&
-\int_0^T\int_0^l [x^*(z,t)qx^*(z,t)+ u^*(z,t)ru^*(z,t)]dzdt\no\\
=&\int_0^Tc\phi(t)[\psi(l,t)+g(l,t)\phi(t)]dt
-\int_0^T\int_0^l [x^*(z,t)\no\\&\times qx^*(z,t)+ u^*(z,t)ru^*(z,t)]dzdt.
\end{align}
Combining (\ref{72}) and (\ref{74}), it derives that
\begin{align}\label{75}
& \int_0^T\int_0^l [x^*(z,t)qx^*(z,t)+ u^*(z,t)ru^*(z,t)]dzdt\no\\
=&\int_0^Tc\phi(t)[\psi(l,t)+g(l,t)\phi(t)]dt\no\\&+\int_0^l \varphi(z)[g(z,0)\varphi(z)+\psi(z,0)]dz\no\\&-\int_0^l x^*(z,T)[p(z)x^*(z,T)+\gamma(z)]dz.
\end{align}

Therefore, by using (\ref{3}), the optimal cost (\ref{2}) is equal to
\begin{align}\label{76}
J^*=& \frac{1}{2}\int_0^T\int_0^l [x^*(z,t)qx^*(z,t)+ u^*(z,t)ru^*(z,t)]dzdt\no\\&+\int_0^l x^*(z,T)p(z)x^*(z,T)dz\no\\
=& \frac{1}{2}\bigg\{\int_0^Tc\phi(t)[\psi(l,t)+g(l,t)\phi(t)]dt\no\\&+\int_0^l \varphi(z)[g(z,0)\varphi(z)+\psi(z,0)]dz\no\\&
- \int_0^l x^*(z,T)\gamma(z)\bigg\}dz,
\end{align}
which implies (\ref{66}).
Thus the proof of Theorem \ref{Thm1} has been completed.

\bibliographystyle{IEEEtran}
\bibliography{2}

\end{document}